\documentclass{amsart}
\usepackage{graphicx}
\newtheorem{thm}{Theorem}[section]

\newtheorem{lem}[thm]{Lemma}
\newtheorem{prop}[thm]{Proposition}
\theoremstyle{definition}

\theoremstyle{remark}

\numberwithin{equation}{section}

\begin{document}

\title[ Submanifolds of a Conformal Sasakian
manifold]
{Submanifolds of a Conformal Sasakian
manifold}
\author{Esmaiel Abedi}%
\address{
Department of Mathematics\\
Azarbaijan shahid Madani University,Tabriz 53751 71379, Iran}%
\email{esabedi@azaruniv.edu}%

%


\date{}

\dedicatory{}

\commby{}


\begin{abstract}
In the present paper, we introduce  conformal Sasakian manifolds. Some results on geometry of such manifolds and their associated submanifolds are provided. 
\end{abstract}

\maketitle

\section{INTRODUCTION}

  A $(2n+1)$-dimensional Riemannian manifold $(M,g)$ said to be a Sasakian manifold if it admits an endomorphism $\varphi$
   of its tangent bundle $TM$, a vector field $\xi$ and a $1$-form $\eta$ satisfying
\begin{eqnarray}
&&\phi^2=-Id+\eta\otimes\xi,\quad  \eta(\xi)=1,\nonumber\\
&& g(\phi X,\phi Y)=g(X,Y)-\eta(X)\eta(Y),\nonumber\\
&&({\nabla}_{X}\phi)Y=g(X,Y)\xi-\eta(Y)X,\nonumber
\end{eqnarray}
for all vector fields $X,Y$ on $M$, where $\nabla$ denotes the Riemannian connection  \cite{Bl}.\\
The close relationship between Kaehler manifolds and Sasakian manifolds naturally leads to the question which objects,
 methods and theorems can be transfered from one to the other.
 Banaru \cite{BT} succeed to classify the sixteen classes of almost Hermitian manifolds by using the two tensors of
Kirichenko, which are called the Kirichenko's tensors. 
The locally conformal Kaehler manifold is one of the sixteen
classes of almost Hermitian manifolds. Libermann did the first study
 on locally conformal Kaehler manifolds \cite{Ph}. Vaisman, put down some geometrical
  conditions for locally conformal Kaehler manifolds \cite{vz}, and Tricerri mentioned different examples about the
  locally conformal Kaehler manifolds \cite{tr}.
The paper is organized as follows:\\
 In section $ 2 $, we  recall som preliminary definitions then we  introduce conformal Sasakian manifolds. Furthermore, we give some basic results on conformal Sasakian manifolds and their submanifolds. In section $ 3 $,   we obtain a  necessary and sufficient condition for the invariant submanifolds of a conformal Sasakian manifold to be minimal.  In  section $ 4 $, we study anti-invariant submanifolds of a conformal Sasakian manifold and obtain the conditions under which these type submanifolds have a flat normal connection. Section $5$ considers  $CR$-submanifolds of a conformal Sasakian manifold with distributions $ D $ and $ D^{\perp} $. we find the conditions under which  $D^{\perp}$ is integrable or totally geodesic.


\section{riemannian geometry of conformal Sasakian manifolds}
A differentiable manifold ${M}^{2n+1}$ is said an almost contact
manifold if it admits a vector field $\xi$, a one-form $\eta$ and
a (1,1)-tensor field $\varphi$ with the
 following properties
\begin{eqnarray}
 \varphi^{2}= -Id + \eta \otimes \xi,\quad\eta(\xi) = 1,\quad\varphi\xi = 0,\quad\eta o\varphi =
 0.
 \end{eqnarray}
 Furthermore, if $M$ be a Riemannian manifold with the Riemannian metric 
  $g$ such that
  \begin{eqnarray}
  g(\varphi X,\varphi Y)={g}(X,Y) - \eta(X)\eta(Y),\nonumber
 \end{eqnarray}
 for all $ X, Y $ on $ M $, then $(\varphi,\xi,\eta,g)$ is called an almost contact metric
 structure on $ M $ and $ M $ is said  an almost contact metric manifold. \\
 A Sasakian manifold is a normal contact metric manifold, that is, an almost contact metric manifold such that $ d\Phi = 0 $ and $ [\varphi ,\varphi](X,Y) = -2d\eta(X,Y)\xi $ for all $ X, Y $ on $ M $, where $ [\varphi ,\varphi] $ is  the Nijenhuis torsion of $ \varphi $.\\
An almost contact metric manifold $({M}^{2n+1},\varphi,\xi,\eta,g)$  is said a
   Sasakian manifold if and only if 
 \begin{eqnarray}
  &&({\nabla}_{X}\varphi)Y= g(X,Y)\xi - \eta(Y)X, 
 \end{eqnarray}
for all vector fields $X, Y$ on $M$, where ${\nabla}$ denotes the Levi-Civita connection
 with respect to $g$ \cite{Bl}.\\
Let $(\acute{M}^{m},\acute{g})$  be a Riemannian manifold into Riemannian manifold $(M^{n},g)$, $ m< n $, with isometric immersion $ \iota:  (\acute{M},\acute{g}) \longrightarrow (M,g)$.
Then the
    Gauss and Weingarten formulas are given by
 \begin{eqnarray}
 &&\nabla_{X}Y= \acute{\nabla}_{X}Y + h(X,Y),\nonumber\\
 &&\nabla_{X}N=-A_{N}X +\nabla^{\perp}_{X}N,\nonumber
 \end{eqnarray}
for all $ X, Y $ tangent to $  \acute{M}$ and normal vector field $ N $ on $ \acute{M} $, where $\acute{\nabla}$ and $\nabla$ are the Levi-Civita connections
of $\acute{M}$ and $M$, respectively, also $h$ and $A_{N}$
are the second fundamental form and the shape operator corresponding to $N$,
respectively and $ \nabla^{\perp} $ is the normal connection on $ T^{\perp}\acute{M} $. 
Let $\acute{R}$ and $R$ denote the curvature tensors on
$\acute{M}$ and $M$, respectively, then
 the Gauss and Codazzi equations are given by
  \begin{eqnarray}
  g(R(X,Y)Z,W)\!\!\!\!&=&\!\!\!\!\acute{g}(\acute{R}(X,Y)Z,W)+h(
  Y,Z)h(X,W)\\
  &&-h(X,Z)h(Y,W),\nonumber\\
 \acute{g}(\acute{R}(X,Y)Z,N)\!\!\!\!&=&\!\!\!\!g((\acute{\nabla}_{X}A_{N})Y -
   (\acute{\nabla}_{Y}A_{N})X,Z)\nonumber\\
&& +\sum_{b=1}^{p}\{S_{ba}(X)g^{\prime}(A_{b}Y,Z)-S_{ba}(Y)g^{\prime}(A_{b}X,Z)\},
\end{eqnarray}
for all $X, Y, Z, W \in T\acute{M}$  and $N_{a} \in T^{\perp}\acute{M}$, where $ A_{a} $ is the shape operator with respect to $ N_{a} $, $ a: 1,\cdots ,p = n-m $.
Also, let $R^{\perp}$ be the normal curvature tensor of $\acute{M}$ then we
will have the Ricci equation by following
\begin{eqnarray}
  g(R(X,Y)N_{1},N_{2})\!\!\!\!&=&\!\!\!\!g(R^{\bot}(X,Y)N_{1},N_{2})-g([A_{1},A_{2}]X,Y),
\end{eqnarray}
where $N_{1}, N_{2}$ are unit normal vector fields on $\acute{M}$ and
$A_{1}, A_{2}$ are the shape operators with respect to $N_{1}, N_{2}$.\\
A smooth manifold ${M}^{2n+1}$ with an almost contact metric structure $ (\varphi,\eta,\xi,g) $ is called a
conformal Sasakian manifold if there is a positive smooth
function $f:M\longrightarrow \mathbb{R}$  such that 
\begin{eqnarray}
\quad\widetilde{g}=exp(f)g,\quad\widetilde{\varphi}=\varphi,\quad\widetilde{\eta}={exp(f)}^{\frac{1}{2}}\eta,
\quad\widetilde{\xi}={exp(-f)}^{\frac{1}{2}}\xi
\end{eqnarray}
is a Sasakian structure on $ M $.
Let $\widetilde{\nabla}$ and $\nabla$
denote connections of $M$  related to metrics
$\widetilde{g}$ and $g$, respectively.
 Using Koszul formula, we
derive the following relation between the connections
 $\widetilde{\nabla}$ and $\nabla$
\begin{eqnarray}
\widetilde{\nabla}_{X}Y=\nabla_{X}Y+\frac{1}{2}\{\omega(X)Y+\omega(Y)X-g(X,Y)\omega^{\sharp}\},
\end{eqnarray}
for all vector fields $ X, Y $ on $ M $, 
so that $\omega(X)=X(f)$ and $\omega^{\sharp}$ is vector field of
metrically equavalente to one form of $\omega$, that is,
$g(\omega^{\sharp},X)=\omega(X)$. Vector field $ \omega^{\sharp} = grad f $ is called the Lee vector field of conformal Sasakian manifold $ M $.\\
Then with a straightforward computation we will have
\begin{eqnarray}
\exp(-f)\widetilde{R}(X,Y,Z,W)&=&R(X,Y,Z,W)\nonumber\\ &+&\frac{1}{2}\{B(X,Z)g(Y,W)-B(Y,Z)g(X,W)\nonumber\\
&+&B(Y,W)g(X,Z)-B(X,W)g(Y,Z)\}\nonumber\\
&+&\frac{1}{4}\|\omega^{\sharp}\|^{2}\{g(X,Z)g(Y,W)-g(Y,Z)g(X,W)\},
\end{eqnarray}
for all vector fields $ X, Y, Z, W $ on $ M $, where $B:=\nabla \omega - \frac{1}{2}\omega\otimes\omega$ and $R,
\widetilde{R}$ are the curvature tensors of $M$ related to
connections of $\nabla$ and $\widetilde{\nabla}$, respectively.
Also, from (2.2) and (2.6) we have
\begin{eqnarray}
(\nabla_{X}\varphi)Y\!\!\!\!&=&\!\!\!\!(\exp(f))^{\frac{1}{2}}\{g(X,Y)\xi-\eta(Y)X\}\\
&&-\frac{1}{2}\{\omega(\varphi Y)X-\omega(Y)\varphi X
+g(X,Y)\varphi\omega^{\sharp}-g(X,\varphi Y)\omega^{\sharp}\},\nonumber\\
\nabla_{X}\xi\!\!\!\!&=&\!\!\!\!(\exp(-f))^{\frac{1}{2}}\varphi
X+\frac{1}{2}\{\eta(X)\omega^{\sharp} -\omega(\xi)X\},
\end{eqnarray}
for all vector fields $  X, Y$ on $ M $. Let $(\acute{M}^{m},\acute{g})$ be a Riemannian manifold into conformal Sasakian manifold
 $ {M}^{2n+1}$ with isometric immersion $ \iota: (\acute{M}^{m},\acute{g})\longrightarrow (M,g)$. Suppose $\acute{\nabla}$ and
$\acute{R}$ are 
the Levi-Civita connection and curvature tensor on $\acute{M}^{m}$, respectively.
We set
 \begin{eqnarray}
&&PX=tan(\varphi X)\hspace{1cm},\hspace{1cm}FX=nor(\varphi X),\nonumber\\
&&tN=tan(\varphi N)\hspace{1cm},\hspace{1cm}fN=nor(\varphi
N),\nonumber
\end{eqnarray}
 for each $X\in T{M^{\prime}}$ and  $N\in TM'^\perp$.
 Then from  (2.9) we get
 \begin{eqnarray}
\nabla_{X}(\varphi Y)\!\!\!\!&=&\!\!\!\!\varphi\nabla_{X}Y+(\exp(f))^{\frac{1}{2}}\{g(X,Y)\xi-\eta(Y)X\}\nonumber\\
&&-\frac{1}{2}\{\omega(\varphi Y)X-\omega(Y)\varphi X
+g(X,Y)\varphi\omega^{\sharp}-g(X,\varphi Y)\omega^{\sharp}\},
\end{eqnarray}
 for all vector field $ X, Y $ on $ \acute{M} $.  Separating the tangential and normal parts from the above equation we will have
  \begin{eqnarray}
  (\acute{\nabla}_{X}P)Y\!\!\!\!&=&\!\!\!\!(\exp(f))^{\frac{1}{2}}\{g(X,Y)\xi^{\top}-\eta(Y)X\}+A_{FY}X+th(X,Y)\nonumber\\
  &&-\frac{1}{2}\{\omega(\varphi Y)X-\omega(Y)PX +g(X,Y)(\varphi \omega^{\sharp})^{\top} - g(X,\varphi Y)\omega^{\sharp^{\top}}\},\\
  (\acute{\nabla}_{X}F)Y\!\!\!\!&=&\!\!\!\!fh(X,Y)-h(X,PY)\nonumber\\
  &&+\frac{1}{2}\{\omega(Y)FX -g(X,Y)(\varphi \omega^{\sharp})^{\bot}+g(X,\varphi Y){\omega^{\sharp}}^{\bot}\},\\
  (\acute{\nabla}_{X}t)N\!\!\!\!&=&\!\!\!\!A_{fN}X-PA_{N}X-\frac{1}{2}\{-\omega(N)PX +\omega(\varphi N)X -g(X,\varphi N)\omega^{\sharp^{\top}}\},\\
  \qquad(\acute{\nabla}_{X}f)N\!\!\!\!&=&\!\!\!\!-h(X,tN)-FA_{N}X+\frac{1}{2}\{\omega(N)FX +g(X,\varphi N){\omega^{\sharp}}^{\bot}\},
 \end{eqnarray}
 for all $X,Y\in T{M^{\prime}}$ and  $N\in TM'^\perp$.\\
 We  need the equations of Gauss, Codazzi and Ricci  between manifolds $\acute{M}^{m}$ and
$ M^{2n+1} $ conformal Sasakian manifold $M$, thus from (2.3), (2.4), (2.5) and (2.8) we get
\begin{eqnarray}
exp(-f)\widetilde{R}(X,Y,Z,W)&=&\acute{R}(X,Y,Z,W)-\frac{1}{2}\{B(X,Z)\acute{g}(Y,W)\nonumber\\
&-&B(Y,Z)\acute{g}(X,W)+B(Y,W)\acute{g}(X,Z)-B(X,W)\acute{g}(Y,Z)\}\nonumber\\
&-&\frac{1}{4}\|\omega^{\sharp}\|^{2}\{\acute{g}(X,Z)\acute{g}(Y,W) - \acute{g}(Y,Z)\acute{g}(X,W)\}\nonumber\\
&+&\sum_{a=1}^{p}\{\acute{g}(A_{a}Y,Z)\acute{g}(A_{a}X,W) - \acute{g}(A_{a}X,Z)\acute{g}(A_{a}Y,W)\},\\
exp(-f)\widetilde{R}(X,Y,Z,N_{a})&=& \acute{g}((\acute{\nabla}_{X}A_{a})Y-(\acute{\nabla}_{Y}A_{a})X,Z)\nonumber\\
&+&\sum_{b=1}^{p}\{S_{ba}(X)\acute{g}(A_{b}Y,Z) - S_{ba}(Y)\acute{g}(A_{b}X,Z)\}\nonumber\\
&+&\frac{1}{2}\{\acute{g}(X,Z)g(\nabla_{Y}\omega^{\sharp},N_{a})-\acute{g}(Y,Z)g(\nabla_{X}\omega^{\sharp},N_{a})\},\\
exp(-f)\widetilde{R}(X,Y,N_{a},N_{b})&=&\acute{g}([A_{2},A_{1}]X,Y)+g(R^{\bot}(X,Y)N_{1},N_{2}),\nonumber
\end{eqnarray}
for all $ X, Y, Z \in T\acute{M} $ and $ N_{a}\in T^{\perp}\acute{M} $, where $ A_{a} $ is the shape operator with respect to $ N_{a} $, $a:1, \cdots ,p=2n-m+1$.
\section{invariant submanifolds}
A submanifold $ \acute{M} $ of a conformal Sasakian manifold $ M $ is called an invariant submanifold of $ M $ if $ \varphi T\acute{M}\subset T\acute{M} $. Hence, $ \varphi N\in T^{\perp}\acute{M} $ for each $ N\in T^{\perp}\acute{M} $, that is, $ tN \equiv 0 $.
\begin{thm}
Let $\acute{M}^{m}$ be an invariant submanifold of a conformal Sasakian
manifold  $M^{2n+1}$ tangent to $\xi$. Then
$\acute{M}$ is minimal if  $\omega^{\sharp}$ 
 is tangent to $\acute{M}$.
\begin{proof}
By relation (2.9) and the Gauss formula we have
\begin{eqnarray}
h(X,\varphi Y)&=&\varphi h(X,Y)-(\acute{\nabla}_{X}\varphi)Y+exp(f)^{\frac{1}{2}}\{g(X,Y)\xi - \eta(Y)X\}\nonumber\\
&-&\frac{1}{2}\{\omega(\varphi Y)X-\omega(Y)\varphi X-g(X,\varphi Y)\omega^{\sharp}+g(X,Y)\varphi\omega^{\sharp}\},
\end{eqnarray}
for all $ X, Y \in T\acute{M} $. Since $ \acute{M} $ is invariant, comparing  tangential  and  normal parts we get
\begin{eqnarray}
&&h(X,\varphi Y)=\varphi
h(X,Y)-\frac{1}{2}\{g(X,Y)(\varphi \omega^{\sharp})^{\perp}-g(X,\varphi Y){\omega^{\sharp}}^{\bot}\}.
\end{eqnarray}
Since $\xi\in T\acute{M}$, taking $ X = \varphi X  $ in (3.2)  we obtain 
\begin{eqnarray*}
&&h(\varphi X,\varphi Y)+h(X,Y)=\{\acute{g}(X,Y)-\frac{1}{2}\eta(X)\eta(Y)\}{\omega^{\sharp}}^{\bot},
\end{eqnarray*}
for all $ X, Y $ on $ \acute{M} $. Again, since $ \xi \in T\acute{M} $, we put $ X = Y = \xi $ in (3.2) then we find
\begin{eqnarray*}
&&h(\xi,\xi)=\frac{1}{2}{\omega^{\sharp}}^{\bot}.
\end{eqnarray*}
Let $\{E_{\alpha},\varphi E_{\alpha}, \xi | \alpha=1,...,n=\frac{m-1}{2}\}$ be an orthonormal frame on
$\acute{M}$ and suppose  $H$ is the mean curvature vector. Then from the above relation we have
\begin{eqnarray*}
H&=&\frac{1}{m}\sum_{\alpha=1}^{n}\{h(\xi,\xi)+h(E_{\alpha},E_{\alpha})+h(\varphi
E_{\alpha},\varphi E_{\alpha})\} \\ 
&=&\frac{1}{m}\{\frac{1}{2}+\sum_{\alpha=1}^{n}\ g(E_{\alpha},E_{\alpha})\}{\omega^{\sharp}}^{\bot}.
\end{eqnarray*}
Thus, theorem  is proved.
\end{proof}
\end{thm}
\section{Anti-invariant submanifolds }
A submanifold $\acute{M}^{m}$ of a conformal Sasakian manifold $ M $ is called an anti-invariant of $ M $ if
$\varphi T\acute{M}\subset T^{\perp}\acute{M}$. Then $ \varphi X\in  T^{\perp}\acute{M}$, for each $ X\in T\acute{M} $, that is, $ P\equiv 0 $.
\begin{lem}
Let $\acute{M}^{m}$ be an $m$-dimensional anti-invariant submanifold of  a conformal
Sasakian manifold ${M}^{2n+1}$. Then
\begin{eqnarray}
A_{\varphi Y}X\!\!\!\!&=&\!\!\!\!-\varphi h(X,Y)-(\exp(f))^{\frac{1}{2}}\{g(X,Y)\xi ^{\top}-\eta(Y)X\}\nonumber\\
&&+\frac{1}{2}\{\omega(\varphi Y)X+g(X, Y)\varphi{\omega^{\sharp}}^{\bot}\},
\end{eqnarray}
and
\begin{eqnarray}
&&\acute{g}([A_{\varphi Z},A_{\varphi W}]X,Y)=g(h(X,W),h(Y,Z))-g(h(Y,W),h(X,Z))\nonumber\\
&-&\frac{1}{2}
\{\acute{g}(Y,Z)\omega(h(X,W))-\acute{g}(Y,W)\omega(h(X,Z))+\acute{g}(X,W)\omega(h(Y,Z))\nonumber\\
&-&\acute{g}(X,Z)\omega(h(Y,W))+\omega(\varphi Z)\Phi(Y,h(X,W))-\omega(\varphi W)\Phi(Y,h(X,Z))\nonumber\\
&+&\omega(\varphi W)\Phi(X,h(Y,W))-\omega(\varphi Z)\Phi(X,h(Y,W))\}\nonumber\\
&-&\frac{1}{4}\{\omega(\varphi W)\omega(\varphi X)\acute{g}(Y,Z)-\omega(\varphi Z)\omega(\varphi X)\acute{g}(Y,W)\nonumber\\
&+&\omega(\varphi Z)\omega(\varphi Y)\acute{g}(X,W)-\omega(\varphi W)\omega(\varphi Y)\acute{g}(X,Z)\nonumber\\
&+&\|{\omega}^{\sharp}\|^{2}\{\acute{g}(X,Z)\acute{g}(Y,W)-\acute{g}(X,W)\acute{g}(Y,Z)\}\}\nonumber\\
&-&\frac{1}{2}(\exp(f))^{\frac{1}{2}}\{2\eta(Z)\Phi(Y,h(X,W))-2\eta(W)\Phi(Y,h(X,Z))\nonumber\\
&+&2\eta(W)\Phi(X,h(Y,Z))-2\eta(Z)\Phi(X,h(Y,W))\nonumber\\
&+&\omega(\varphi Z)\eta(Y)\acute{g}(X,W)-\omega(\varphi W)\eta(Y)g^{\prime}(X,Z)\nonumber\\
&+&\omega(\varphi X)\eta(W)\acute{g}(Y,Z)-\omega(\varphi X)\eta(Z)\acute{g}(Y,W)\nonumber\\
&+&\omega(\varphi W)\eta(X)\acute{g}(Y,Z)-\omega(\varphi Z)\eta(X)\acute{g}(Y,W)+\omega(\varphi Y)\eta(Z)\acute{g}(X,W)\nonumber\\
&-&\omega(\varphi Y)\eta(W)\acute{g}(X,Z)\}\nonumber\\
&+&\exp(f)\{\acute{g}(Y,Z)\acute{g}(X,W)-\acute{g}(X,Z)\acute{g}(Y,W)+\acute{g}(X,Z)\eta(Y)\eta(W)\nonumber\\
&-&\acute{g}(X,W)\eta(Y)\eta(Z)+\acute{g}(Y,W)\eta(X)\eta(Z)-\acute{g}(Y,Z)\eta(X)\eta(W)\},
\end{eqnarray}
for all $X, Y, Z, W\in  T\acute{M}$, where $ \Phi(X,Y)=g(X,\varphi Y)$.
\begin{proof}
Since $P \equiv 0$ then (4.1) follows from (2.12), easily. Also, substituting (4.1) in $ \acute{g}([A_{\varphi Z},A_{\varphi W}]X,Y) = \acute{g}(A_{\varphi W}X,A_{\varphi Z}Y) - \acute{g}(A_{\varphi Z}X,A_{\varphi W}Y) $, we get  (4.2).
\end{proof}
\end{lem}
\begin{prop}\label{prop.1}
Let $\acute{M}^{m}$ be an  anti-invariant
submanifold of  a conformal Sasakian manifold
${M}^{2n+1}$ tangent to $ \xi $. Then $\acute{M}$ has a
flat normal connection if and only if
\begin{eqnarray}
\acute{R}(X,Y)Z&=&\eta(R(X,Y)Z)\xi\nonumber\\
&+&\frac{1}{2}
\{\acute{B}(Y,Z)X-\acute{B}(X,Z)Y+\acute{g}(Y,Z)\acute{B}(X,.)^{\sharp}-\acute{g}(X,Z)\acute{B}(Y,.)^{\sharp}\nonumber\\
&+&B(X,Z)\eta(Y)\xi -B(Y,Z)\eta(X)\xi\nonumber\\
&+&B(Y,\xi)\acute{g}(X,Z)\xi -B(X,\xi)\acute{g}(Y,Z)\xi\}\nonumber\\
&+&\frac{1}{4}\|{\acute{\omega}}^{\sharp}\|^{2}\{\acute{g}(Y,Z)X-\acute{g}(X,Z)Y\}\nonumber\\
&+&\frac{1}{4}\|{\omega}^{\sharp}\|^{2}\{\acute{g}(X,Z)\eta(Y)\xi -\acute{g}(Y,Z)\eta(X)\xi\}\nonumber\\
&+&\{\acute{g}(Y,Z)X-\acute{g}(X,Z)Y+\acute{g}(X,Z)\eta(Y)\xi -\acute{g}(Y,Z)\eta(X)\xi\}\nonumber\\
&-&\exp(f)\{\acute{g}(Y,Z)X-\acute{g}(X,Z)Y\}
\end{eqnarray}
for all $X, Y, Z\in  T\acute{M}$, where $\acute{\omega}^{\sharp}= \omega^{\sharp^{\top}}$ and $\acute{B}=B+\omega o h$.
\begin{proof}
Since $(\tilde{\nabla}_{X}\varphi)Y= \widetilde{g}(X,Y)\xi-\eta(Y) X $ it follows that \cite{Bl}
\begin{eqnarray}
\widetilde{R}(X,Y)\varphi Z\!\!\!\!&=&\!\!\!\!\varphi \widetilde{R}(X,Y)Z-\widetilde{g}(Y,Z)\varphi X+\widetilde{g}(X,Z)\varphi Y\nonumber\\
&&-\widetilde{g}(\varphi Y,Z)X+\widetilde{g}(\varphi X,Z)Y,
\end{eqnarray}
for all $ X, Y, Z \in T\acute{M} $.
Replacing (2.8) in (4.4) we can write
\begin{eqnarray}
R(X,Y)\varphi Z&=&\varphi R(X,Y)Z\nonumber\\
&-&\frac{1}{2}
\{B(X,\varphi Z)Y-B(Y,\varphi Z)X+B(Y,Z)\varphi X\nonumber\\
&-&B(X,Z)\varphi Y+B(Y,.)^{\sharp}g(X,\varphi Z)\nonumber\\
&-&B(X,.)^{\sharp}g(Y,\varphi Z)-\varphi B(Y,.)^{\sharp}g(X,Z)+\varphi B(X,.)^{\sharp}g(Y,Z)\}\nonumber\\
&-&(\frac{1}{4}\|{\omega}^{\sharp}\|^{2}+1)\{g(Y,Z)\varphi X-g(X,Z)\varphi Y+g(X,\varphi Z) X-g(Y,\varphi Z)X\},
\end{eqnarray}
for all  $X, Y, Z, W\in TM^{\prime}$, where $ B(X,Y) = g(B(X,.)^{\sharp},Y) $. Taking the inner product from (4.5) with $\varphi W$
and using the Ricci and Gauss equations, we obtain
\begin{eqnarray}
&&g(R^{\bot}(X,Y)\varphi Z,\varphi W)-\acute{g}([A_{\varphi Z},A_{\varphi W}]X,Y)\nonumber\\
&=&\acute{g}(\acute{R}(X,Y)Z,W)-g(h(X,W),h(Y,Z))\nonumber\\
&+&g(h(Y,W),h(X,Z))-\eta(R(X,Y)Z)\eta(W)\nonumber\\
&-&\frac{1}{2}
\{B(Y,Z)\acute{g}(X,W)-B(X,Z)\acute{g}(Y,W)\nonumber\\
&+&B(X,W)\acute{g}(Y,Z)-B(Y,W)\acute{g}(X,Z)+B(X,Z)\eta(Y)\eta(W)\nonumber\\
&-&B(Y,Z)\eta(X)\eta(W)+ B(Y,\xi)g(X,Z)\eta(W)- B(X,\xi)g(Y,Z)\eta(W)\}\nonumber\\
&-&(\frac{1}{4}\|{\omega}^{\sharp}\|^{2}+1)\{\acute{g}(Y,Z)\acute{g}(X,W)-\acute{g}(X,Z)\acute{g}(Y,W)\nonumber\\
&+&\acute{g}(X,Z)\eta(Y)\eta(W)-\acute{g}(Y,Z)\eta(X)\eta(W)\},
\end{eqnarray}
for all  $X, Y, Z, W\in T\acute{M}$. From (4.1) we get
\begin{eqnarray}
\Phi(Y,h(X,Z))\!\!\!\!&=&\!\!\!\!\Phi(Z,h(X,Y))-(\exp(f))^{\frac{1}{2}}\{\acute{g}(X,Z)\eta(Y)-\acute{g}(X,Y)\eta(Z)\}\nonumber\\
&&+\frac{1}{2}\{\omega(\varphi Z)\acute{g}(X,Y)-\omega(\varphi Y)\acute{g}(X,Z)\},
\end{eqnarray}
for all  $X, Y, Z, W\in T\acute{M}$. Putting  (4.2) into (4.6) and using (4.7), we find
\begin{eqnarray}
-\varphi R^{\bot}(X,Y)\varphi Z\!\!\!\!&=&\!\!\!\! \acute{R}(X,Y)Z-\eta(R(X,Y)Z)\xi\nonumber\\
&-&\frac{1}{2}
\{\acute{B}(Y,Z)X-\acute{B}(X,Z)Y+\acute{g}(Y,Z)\acute{B}(X,.)^{\sharp}\nonumber\\
&-&\acute{g}(X,Z)\acute{B}(Y,.)^{\sharp}+B(X,Z)\eta(Y)\xi -B(Y,Z)\eta(X)\xi\nonumber\\
&+&B(Y,\xi)\acute{g}(X,Z)\xi -B(X,\xi)\acute{g}(Y,Z)\xi \}\nonumber\\
&-&\frac{1}{4}\| \acute{\omega}^{\sharp}\|^{2}\{\acute{g}(Y,Z)X-\acute{g}(X,Z)Y\}\nonumber\\
&-&\frac{1}{4}\| \omega^{\sharp}\|^{2}\{\acute{g}(X,Z)\eta(Y)\xi-\acute{g}(Y,Z)\eta(X)\xi \}\nonumber\\
&-&\{\acute{g}(Y,Z)X-\acute{g}(X,Z)Y+\acute{g}(X,Z)\eta(Y)\xi -\acute{g}(Y,Z)\eta(X)\xi\}\nonumber\\
&+&\exp(f)\{\acute{g}(Y,Z)X-\acute{g}(X,Z)Y\},
\end{eqnarray}
for all  $X, Y, Z, W\in T\acute{M}$. Thus $R^{\bot}=0$ if and only if (4.3) holds.
\end{proof}
\end{prop}
Let ${\acute{M}}^{m}$ be an  anti-invariant
submanifold of a  conformal Sasakian manifold
$M^{2n+1}$. The normal curvature tensor $ R^{\perp} $ of  $\acute{M}$ is called reccurent if
\begin{eqnarray}
R^{\bot}(X,Y)N=\theta(X,Y)N,
\end{eqnarray}
for all $X, Y\in T\acute{M}$and  $N\in T\acute{M}^{\perp}$ holds on $ \acute{M} $, where $\theta$ is a $ 2$-form  on $\acute{M}$.
\begin{thm}\label{thm.1}
Let $\acute{M}^{m}$ be an  anti-invariant
submanifold of a  conformal Sasakian manifold
$M^{2n+1}$ normal to $ \xi $ with reccurent nomal curvature tensor. Then $\acute{M}$ has a flat normal connection.
\begin{proof}
Since $ R^{\perp} $ is reccurent, by (4.9) and using (4.3) in Proposition $4.2$   we obtain
\begin{eqnarray}
\acute{R}(X,Y)Z\!\!\!\!&=&\!\!\!\!\theta(X,Y)Z-\theta(X,Y)\eta(Z)\xi+\eta(R(X,Y)Z)\xi\nonumber\\
&&+\frac{1}{2}
\{\acute{B}(Y,Z)X-\acute{B}(X,Z)Y+\acute{g}(Y,Z)\acute{B}(X,.)^{\sharp}\nonumber\\
&&-\acute{g}(X,Z)\acute{B}(Y,.)^{\sharp}+B(X,Z)\eta(Y)\xi-B(Y,Z)\eta(X)\xi\nonumber\\
&&+B(Y,\xi)\acute{g}(X,Z)\xi-B(X,\xi)\acute{g}(Y,Z)\xi\}\nonumber\\
&&+\frac{1}{4}\|{\acute{\omega}}^{\sharp}\|^{2}\{\acute{g}(Y,Z)X-\acute{g}(X,Z)Y\}\nonumber\\
&&+\frac{1}{4}\|{\omega}^{\sharp}\|^{2}\{\acute{g}(X,Z)\eta(Y)\xi-\acute{g}(Y,Z)\eta(X)\xi\}\nonumber\\
&&+\{\acute{g}(Y,Z)X-\acute{g}(X,Z)Y+\acute{g}(X,Z)\eta(Y)\xi -\acute{g}(Y,Z)\eta(X)\xi\}\nonumber\\
&&-\exp(f)\{\acute{g}(Y,Z)X-\acute{g}(X,Z)Y\},
\end{eqnarray}
for all $ X, Y, Z\in T\acute{M} $. Since $ \xi \in T^{\perp}\acute{M} $, taking the inner product from the above equation with each vector field $ W\in T\acute{M} $ and  Contracting it over $ Z $ and $ W $ we get
\begin{eqnarray}
m\theta(X,Y) = 0,
\end{eqnarray}
for all $ X, Y $ on $ \acute{M} $. Then (4.9) results $R^{\perp}=0$. Thus, $ \acute{M} $ has a flat normal connection.
\end{proof}
\end{thm}
\section{Distribution on submanifolds}
Let ${M}^{2n+1}$ be a conformal Sasakian manifold.
Then $\acute{M}^{m}$ is said a $CR$-submanifold in $ M $ if there exist two
orthogonal complementray distributions $D$ and $D^{\bot}$ of
$T\acute{M}$ such that $ \xi \in T\acute{M}$ and \\
(1) $D$ is invariant by $\varphi$, i.e. $\varphi(D_{p})\subset
D_{p}, \forall p\in \acute{M}$.\\
(2) $D^{\bot}$ is anti-invariant by $\varphi$, i.e.
$\varphi(D_{p}^{\perp})\subset T_{p}^{\bot}\acute{M}  \forall p\in
\acute{M}$.
\begin{thm}\label{thm.1}
Let $(\acute{M}^{m},D)$ be a $CR$-submanifold of a conformal
Sasakian manifold  ${M}^{2n+1}$. Then the
anti-invariant distribution $D^{\bot}$ of $M^{\prime}$ is
integrable.
\begin{proof}
Since $ \Phi(X,Y) = g(X,\varphi Y) $ for all $ X, Y \in T\acute{M}$, we get $\Phi(X,Y)=0$ and
$\Phi(Z,W)=0$ for all $X\in D$ and $Y, Z\in D^{\bot}$.  Since $(M^{2n+1},\varphi, \widetilde{\xi}, \widetilde{\eta}, \widetilde{g})$ is a  Sasakian manifold, we have $d\widetilde{\Phi}=0$, where $\widetilde{\Phi}(X,Y)=\widetilde{g}(X,\varphi Y)$. Thus, we find
\begin{eqnarray}
0&=&d\widetilde{\Phi}\nonumber\\
&=&\!\!\!\!d(\exp(f))\wedge\Phi+\exp(f)d\Phi\nonumber\\
&=&\!\!\!\!\exp(f)(\omega \wedge\Phi+d\Phi).\nonumber
\end{eqnarray}
Using $(\Phi\wedge\omega)(X,Y,Z)=0$ for all $X\in D$ and $Y, Z\in D^{\bot}$ in the above equation we can write
\begin{eqnarray}
0\!\!\!\!&=&\!\!\!\!3(d\Phi)(X,Y,Z)\nonumber\\
&=&\!\!\!\!X(\Phi(Y,Z))+Z(\Phi(Z,X))+W(\Phi(X,Y))\nonumber\\
&&-\Phi([X,Y],Z)-\Phi([Z,X],Y)-\Phi([Y,Z],X)\nonumber\\
&=&\!\!\!\!-g([Y,Z],\varphi X),\nonumber
\end{eqnarray}
hence, $[Y,Z]\in D^{\bot}$ for all $ Y, Z\in D^{\perp} $.
\end{proof}
\end{thm}
Let $(\acute{M}^{m},D)$ be a $CR$-submanifold of a conformal Sasakian manifold
$M^{2n+1}$. Then
 $M^{\prime}$ is said to be mixed totally geodesic if $th(X,Y)=0$
for each $X\in D$ and  $Y\in D^{\bot}$.\cite{Dr}
\begin{thm}\label{thm.1}
Let $({\acute{M}}^{m},D)$ be a $CR$-submanifold of a conformal
Sasakian manifold $M^{2n+1}$ normal to $ \omega^{\sharp} $. Then
$\acute{M}$ is mixed totally geodesic if and only if each
leaf of the anti-invariant distribution $D^{\bot}$ is a  totally
geodesic submanifold  of $ \acute{M}$
\begin{proof}
Let $S$ be a leaf of $D^{\bot}$.   Making use of  the Gauss formula we have
\begin{eqnarray}
h_{S}(Y,Z)=(\acute{\nabla}_{Y}Z)_{D}
\end{eqnarray}
for all $Y, Z\in T(S)=D^{\bot}$ and $X\in D$, where $h_{S}$  is the second fundamental form of $S$ in $\acute{M}$. 
Hence we have
\begin{eqnarray}
-\acute{\nabla}((\acute{\nabla}_{Y}Z)_{D},\varphi
X)=\acute{g}(th(X,Y),Z)-\frac{1}{2}\acute{g}(Y,Z)\omega(\varphi
X ).
\end{eqnarray}
Since $\omega^{\sharp}$ is normal to $ \acute{M} $,in view of (5.2) we get
\begin{eqnarray}
g^{\prime}((\nabla^{\prime}_{Y}Z)_{D},\varphi
X)&=&-g^{\prime}(th(X,Y),Z),
\end{eqnarray}
for all $Y, Z\in T(S) = D^{\bot}$ and $X\in D$.
 So (5.1) and (5.3) complete the proof of the theorem.
\end{proof}
\end{thm}



\end{document}